\documentclass[a4paper]{amsart}

\usepackage[l2tabu,orthodox,abort]{nag}

\usepackage[british]{babel}

\usepackage{mlmodern}
\DeclareFontFamily{OMX}{mlmex}{}
\DeclareFontShape{OMX}{mlmex}{m}{n}{<->mlmex10}{}
\usepackage[T1]{fontenc}
\usepackage[babel,stretch=10,shrink=10,verbose=errors,selected]{microtype}
\usepackage[strict,english=british]{csquotes}

\usepackage{float}
\usepackage{multirow}

\usepackage{amsmath,amssymb,amsthm}
\usepackage{mathtools}
\usepackage[all]{onlyamsmath}
\usepackage{mleftright}

\usepackage[
	backend=biber,
	style=numeric-comp,
	sorting=nyt,
	sortcites=true,
	maxnames=5,
	giveninits=true,
	date=year,
	useprefix=true
]{biblatex}
\DeclareFieldFormat{eid}{no.~#1}
\usepackage{enumitem}

\usepackage[nobiblatex]{xurl}
\usepackage[breaklinks=true,pdflang=en-GB,colorlinks=true]{hyperref}
\usepackage{xcolor}
\colorlet{citecolor}{green!75!black}
\colorlet{linkcolor}{red!75!black}
\colorlet{urlcolor}{blue!75!black}
\hypersetup{citecolor=citecolor,linkcolor=linkcolor,urlcolor=urlcolor}
\usepackage{zref-clever}
\zcsetup{cap,nameinlink=false}
\usepackage{keytheorems}

\usepackage{orcidlink}
\usepackage{ellipsis}

\numberwithin{equation}{section}

\theoremstyle{plain}
\newtheorem{theorem}{Theorem}[section]
\newtheorem{lemma}[theorem]{Lemma}
\newtheorem{corollary}[theorem]{Corollary}
\newtheorem{proposition}[theorem]{Proposition}
\newtheorem{conjecture}[theorem]{Conjecture}

\theoremstyle{definition}
\newtheorem{definition}[theorem]{Definition}

\newtheorem{question}[theorem]{Question}

\theoremstyle{remark}

\DeclareMathOperator{\End}{End}
\DeclareMathOperator{\Aut}{Aut}
\DeclareMathOperator{\Inn}{Inn}
\DeclareMathOperator{\Out}{Out}
\DeclareMathOperator{\Hol}{Hol}
\newcommand{\Autc}{\Aut_{\mathrm{c}}}
\newcommand{\Outc}{\Out_{\mathrm{c}}}

\DeclareMathOperator{\im}{im}

\DeclareMathOperator{\id}{id}

\newcommand{\R}{\mathfrak{R}}
\DeclareMathOperator{\Fix}{Fix}
\DeclareMathOperator{\Spec}{Spec}
\DeclareMathOperator{\ESpec}{ESpec}
\newcommand{\SpecR}[1]{\Spec_R(#1)}
\newcommand{\ESpecR}[1]{\ESpec_R(#1)}

\newcommand{\NN}{\mathbb{N}}
\newcommand{\ZZ}{\mathbb{Z}}

\newcommand{\bvarphi}{{\bar{\varphi}}}

\newcommand{\card}[1]{\# {#1}} 
\newcommand{\SmallGroup}[2]{\texttt{[#1,#2]}}
\newcommand{\assign}{\vcentcolon=}

\title{Extreme Reidemeister spectra of finite groups}
\author[Sam Tertooy]{Sam Tertooy\,\orcidlink{0000-0002-5750-9153}}
\date{\today}
\address{KU Leuven, Kulak Kortrijk Campus\\
	E.~Sabbelaan 53\\
	8500 Kortrijk\\
	Belgium}
\email{\href{mailto:sam.tertooy@kuleuven.be}{sam.tertooy@kuleuven.be}}
\urladdr{\url{https://stertooy.github.io}}

\subjclass[2020]{Primary: 20D45; Secondary: 20E45}
\keywords{Twisted conjugacy, finite groups, Reidemeister number, Reidemeister spectrum}

\begin{document}
	
\begin{abstract}
We extend the notions of ``\(R_\infty\)-property'' and ``full (extended) Reidemeister spectrum'' to finite groups in a meaningful way. We provide examples of finite groups admitting these properties, if they exist, by looking at groups of small order as well as (quasi)simple groups.
\end{abstract}
	
\maketitle

	\begin{center}
	This is an Accepted Manuscript of an article published by Elsevier in Topology and its Applications on 1 Jan 2025, available online:  \href{https://doi.org/nv2n}{https://doi.org/nv2n}.
	\end{center}

\section{Introduction}
Let \(G\) be a group and let \(\varphi\colon G\to G\) be an endomorphism. Two elements \(g_1,g_2\) of \(G\) are said to be \emph{\(\varphi\)-twisted conjugate} if there exists some \(h \in G\) such that \(g_1 = h g_2 \varphi(h)^{-1}\). This creates an equivalence relation, the classes of which are called \(\varphi\)-twisted conjugacy classes or Reidemeister classes; the Reidemeister class of \(g \in G\) is denoted by \([g]_{\varphi}\). We write \(\R(\varphi)\) for the set of all \(\varphi\)-twisted conjugacy classes. The cardinality of this set is called the Reidemeister number and is denoted by \(R(\varphi)\).

If \(\varphi = \id_G\), then \(\varphi\)-twisted conjugation is just ordinary conjugation. In that case, we will denote the conjugacy class of \(g \in G\) by \([g]\), the set of all conjugacy classes of \(G\) by \(\mathcal{C}(G)\) and the number of conjugacy classes by \(k(G)\), also known as the class number of \(G\).

A common objective in the study of twisted conjugacy is to determine the \emph{Reidemeister spectrum} \(\SpecR{G}\) of a group \(G\), i.e. the set of all possible Reidemeister numbers of automorphisms of \(G\). Two extreme cases can be discerned. The first case is that all Reidemeister numbers of automorphisms are infinite. In that case, we say \(G\) has the \emph{\(R_\infty\)-property}, a term coined by Taback and Wong in \cite{tw07-a}. The second case is that every positive integer, as well as infinity, appears as the Reidemeister number of some automorphism of the group. In that case we say \(G\) has \emph{full} Reidemeister spectrum. One can also consider the set of Reidemeister numbers of all endomorphisms rather than just automorphisms, which is called the \emph{extended Reidemeister spectrum} \(\ESpecR{G}\). Similar extreme cases can be studied here, i.e. when this spectrum is either \(\{1,\infty\}\) or \(\NN \cup \{\infty\}\).

As the presence of infinity in these definitions suggests, the study of twisted conjugacy on groups has been strongly focused on infinite groups --- none of the extreme cases can happen for a finite group. While the first results for twisted conjugacy on finite groups, by Fel'shtyn and Hill, date back to the 90's \cite{fh94-a}, it wasn't until very recently that finite groups appeared as the central object of study, particularly in the work of Senden \cite{send21-a,send23-a,send23-b}.

In the present manuscript, we define notions for the extreme cases for \(\SpecR{G}\) and \(\ESpecR{G}\) if \(G\) is finite. With the help of GAP \cite{gap25-a} and the SmallGrp \cite{beo24-a}, SmallClassNr \cite{tert25-c} and TwistedConjugacy \cite{tert26-a} packages, we provide some preliminary results, and finally we propose some questions that would further the study of Reidemeister spectra of finite groups.

\section{Preliminaries}
\label{sec:prelims}

\begin{definition}
	An endomorphism \(\varphi\) of a group \(G\) is called \emph{fixed-point-free} if it does not preserve any non-trivial element of \(G\), i.e. \(\varphi(g) \neq g\) for all \(g \neq 1\).
\end{definition}

The stabiliser of the \(\varphi\)-twisted conjugacy action of a group \(G\) on the identity \(1\) is exactly \(\Fix(\varphi)\), the set of elements of \(g\) for which \(\varphi(g) = g\). If \(G\) is finite, we can apply the orbit-stabiliser theorem to find
\begin{equation*}
	\card{[1]_\varphi} \cdot \card{\Fix(\varphi)} = \card{G},
\end{equation*}
which leads to the following nice result.
\begin{proposition}
	Let \(G\) be a finite group and let \(\varphi \in \End(G)\). Then
	\begin{equation*}
		\varphi \text{ is fixed-point-free} \iff R(\varphi) = 1.
	\end{equation*}
\end{proposition}

\begin{definition}
	An endomorphism \(\varphi\) of a group \(G\) is called \emph{class-preserving} if it preserves the conjugacy classes of \(G\), i.e. \([\varphi(g)] = [g]\) for all \(g \in G\).
\end{definition}
A class-preserving endomorphism \(\varphi\) must be injective. Therefore, if \(G\) is finite, \(\varphi\) will actually be an automorphism. The set of class-preserving automorphisms of a group \(G\) is denoted by \(\Autc(G)\) and forms a normal subgroup of \(\Aut(G)\) containing \(\Inn(G)\); the quotient \(\Autc(G)/\Inn(G)\) is denoted by \(\Outc(G)\).

\medskip

An endomorphism \(\varphi\) of a group \(G\) induces a self-map \(\Phi\) on the set \(\mathcal{C}(G)\) of conjugacy classes of \(G\), given by
\begin{equation*}
	\Phi\colon \mathcal{C}(G) \to \mathcal{C}(G)\colon [g] \mapsto [\varphi(g)].
\end{equation*}
The connection between this \(\Phi\) and the Reidemeister number \(R(\varphi)\) was established by Fel'shtyn and Hill in \cite[Thm.~5]{fh94-a}.

\begin{theorem}
	\label{thm:reidnrfixedpoints}
	Let \(G\) be a finite group, let \(\varphi \in \End(G)\) and let \(\Phi\) be the self-map on \(\mathcal{C}(G)\) induced by \(\varphi\). Then
	\begin{equation*}
		R(\varphi) = \card{\Fix(\Phi)}.
	\end{equation*}
\end{theorem}

We thus obtain a natural boundary on Reidemeister number for finite groups.
\begin{corollary}
	\label{cor:RiskGiffclasspreserving}
	Let \(G\) be a finite group with class number \(k(G)\) and let \(\varphi \in \End(G)\). Then
	\begin{equation*}
		R(\varphi) \leq k(G),
	\end{equation*}
	with equality holding if and only if \(\varphi\) is class-preserving.
\end{corollary}

Another consequence of \zcref{thm:reidnrfixedpoints} is that Reidemeister numbers of odd order groups must be odd, a result that was also obtained in \cite[Prop.~8.2.3]{send23-a}.
\begin{corollary}
	\label{cor:oddordergroupreidnrodd}
	Let \(G\) be a finite group of odd order, and let \(\varphi \in \End(G)\). Then \(R(\varphi)\) is odd.
\end{corollary}
\begin{proof}
	Let \(g \in G\) be non-trivial. Then \(g\) must have odd order and hence \(g \neq g^{-1}\). Now suppose, by contradiction, that \([g] = [g^{-1}]\). Due to the orbit-stabiliser theorem, \([g]\) contains an odd number of elements. But for any \(h \in [g]\) we must have that \(h^{-1} \in [g^{-1}] = [g]\), and therefore \([g]\) contains an even number of elements. Thus, \([g]\) and \([g^{-1}]\) are distinct conjugacy classes.
	
	If \(\varphi \in \End(G)\), then \([\varphi(g)] = [g]\) if and only if \([\varphi(g^{-1})] = [g^{-1}]\). So \(\varphi\) fixes the conjugacy class \([1]\), as well as a certain number of pairs of conjugacy classes \([g],[g^{-1}]\), and therefore \(R(\varphi)\) is odd.
\end{proof}

The following result can be extracted from \cite[Sec.~4.2]{fz20-b}.
\begin{proposition}
	\label{prop:goodnormalsubgroup}
	Let \(G\) be a group and \(\varphi\) an endomorphism of \(G\). Then the set
	\begin{equation*}
		N_\varphi \assign \bigcup_{n \in \NN} \ker(\varphi^n) = \{g \in G \mid \exists n \in \NN : \varphi^n(g) = 1\}
	\end{equation*}
	is a \(\varphi\)-invariant normal subgroup of \(G\). Moreover, if \(\bvarphi\) denotes the induced endomorphism on \(G/N_\varphi\), then \(\bvarphi\) is injective and \(R(\varphi) = R(\bvarphi)\).
\end{proposition}
We remark that if \(G\) is finite in the above \zcref[noref,nocap]{prop:goodnormalsubgroup}, then \(\bvarphi\) is an automorphism.

\section{Trivial Reidemeister spectrum}

A group \(G\) is said to have the \(R_\infty\)-property if \(\SpecR{G} = \{\infty\}\), which obviously can never happen for a finite group. For automorphisms of infinite groups, infinity is a natural upper bound for Reidemeister numbers, and for many groups this bound is attained by the identity automorphism.

For a finite group \(G\), \zcref{cor:RiskGiffclasspreserving} shows us that the class number \(k(G)\) fulfils the role of natural upper bound, and again it is attained by the identity automorphism. This lets us define an analogue to the \(R_\infty\)-property for finite groups.

\begin{definition}
	A finite group \(G\) is said to have \emph{trivial Reidemeister spectrum} if
	\begin{equation*}
		\SpecR{G} = \{k(G)\}.
	\end{equation*}
\end{definition}

We can quickly come up with the following characterisation for finite groups with trivial Reidemeister spectrum.
\begin{proposition}
	A finite group \(G\) has trivial Reidemeister spectrum if and only if every automorphism of \(G\) is class-preserving.
\end{proposition}
In \cite{mann99-a}, Mann asked the following question: ``Do all \(p\)-groups have automorphisms that are not class-preserving? If the answer is no, which are the groups that have only class-preserving automorphisms?'' When this question was asked, there were already examples in the literature showing that the first part has a negative answer, e.g. \cite{hein79-a,mali92-a}. The second part of this question, when applied to all finite groups and not just \(p\)-groups, can be reformulated as follows.
\begin{question}
	Which finite groups have trivial Reidemeister spectrum?
\end{question}

The easiest examples are groups with trivial outer automorphism group, e.g. the symmetric groups \(S_n\) (\(n \neq 6\)) and the holomorphs \(A \rtimes \Aut(A)\) with \(A\) abelian of odd order \cite[\S 1]{mill08-a}. A theorem of Wielandt \cite[Satz (45)]{wiel39-a} says that for any centreless group \(G\), if we define the sequence \(G_1 \assign G\), \(G_{i+1} \assign \Aut(G_i)\), then this sequence eventually becomes constant at some \(G_n\). This \(G_n\) is then a \emph{complete group}, i.e. \(Z(G_n) = 1\) and \(\Out(G_n) = 1\), and hence has trivial Reidemeister spectrum.

For finite abelian groups, the Reidemeister spectra were completely determined by Senden in \cite{send23-b}; the only non-trivial abelian group with trivial Reidemeister spectrum is \(\ZZ_2\). A result of Gaschütz \cite{gasc66-a} says that non-abelian \(p\)-groups have non-trivial outer automorphism group. Thus, the only way for a \(p\)-group \(G\) to have trivial Reidemeister spectrum, is if \(\Inn(G) \subsetneq \Autc(G) = \Aut(G)\). We list some examples of such \(p\)-groups, where \SmallGroup{n}{k} denotes the ID of the group in the SmallGrp library. In particular, \texttt{n} is the order of the group.
\begin{itemize}
	\item \SmallGroup{128}{932}, a 2-group of class \(4\), found in \cite[Sec.~4.3]{mizu12-a}.
	\item For all \(p \geq 3\), there exists such \(p\)-group of class \(3\) and order \(p^6\). The groups of order \(3^6\) are \SmallGroup{729}{224}, \SmallGroup{729}{232} and \SmallGroup{729}{234} and were obtained in \cite[Sec.~4.3]{mizu12-a}; for \(p \geq 5\) this result was obtained in \cite{mali92-a}.
	\item For all \(p \geq 3\) and \(n \in \NN\), there exists such \(p\)-group of class \(2\) and order \(p^{6n+3}\) if \(p > 3\), or order \(p^{6n+9}\) if \(p = 3\) \cite{hein79-a}.
\end{itemize}
The list of all finite groups of order \(< 512\) with trivial Reidemeister spectrum can be found in \zcref{tbl:trivreidspecloworder}.

\begin{table}[!ht]
	\centering
	\small
	\caption{Groups of order \(< 512\) with trivial Reidemeister spectrum}
	\label{tbl:trivreidspecloworder}
	\begin{tabular}{c|c|c|r}
       SmallGrp ID & Structure & \(\card{\Out(G)}\) & \(k(G)\) \\ \hline
\texttt{[ \ \ 1, \ \ \ 1 ]} & \(1\) & 1 & 1 \\ 
\texttt{[ \ \ 2, \ \ \ 1 ]} & \(\ZZ_2\) & 1 & 2 \\ 
\texttt{[ \ \ 6, \ \ \ 1 ]} & \(S_3 \cong \Hol(\ZZ_3)\) & 1 & 3 \\ 
\texttt{[ \ 20, \ \ \ 3 ]} & \(\Hol(\ZZ_5)\) & 1 & 5 \\
\texttt{[ \ 24, \ \ 12 ]} & \(S_4\) & 1 & 5 \\ 
\texttt{[ \ 42, \ \ \ 1 ]} & \(\Hol(\ZZ_7)\) & 1 & 7 \\ 
\texttt{[ \ 54, \ \ \ 6 ]} & \(\Hol(\ZZ_9)\) & 1 & 10 \\
\texttt{[ 110, \ \ \ 1 ]} & \(\Hol(\ZZ_{11})\) & 1 & 11 \\
\texttt{[ 120, \ \ 34 ]} & \(S_5\) & 1 & 7 \\ 
\texttt{[ 120, \ \ 36 ]} & \(\Hol(\ZZ_{15})\) & 1 & 15 \\ 
\texttt{[ 128, \ 932 ]} &  & \textbf{4} & 17 \\ 
\texttt{[ 144, \ 182 ]} &  & 1 & 9 \\ 
\texttt{[ 144, \ 183 ]} & \(S_3 \times S_4\) & 1 & 15 \\ 
\texttt{[ 156, \ \ \ 7 ]} & \(\Hol(\ZZ_{13})\) & 1 & 13 \\ 
\texttt{[ 168, \ \ 43 ]} &  & 1 & 8 \\ 
\texttt{[ 216, \ \ 90 ]} &  & 1 & 19 \\ 
\texttt{[ 252, \ \ 26 ]} & \(\Hol(\ZZ_{21})\) & 1 & 21 \\ 
\texttt{[ 272, \ \ 50 ]} & \(\Hol(\ZZ_{17})\) & 1 & 17 \\ 
\texttt{[ 320, 1635 ]} &  & 1 & 11 \\ 
\texttt{[ 324, \ 118 ]} & \(S_3 \times \Hol(\ZZ_9)\) & 1 & 30 \\
\texttt{[ 336, \ 208 ]} &  & 1 & 9 \\ 
\texttt{[ 336, \ 210 ]} &  & 1 & 16 \\ 
\texttt{[ 342, \ \ \ 7 ]} & \(\Hol(\ZZ_{19})\) & 1 & 19 \\ 
\texttt{[ 384, 5677 ]} &  & 1 & 16 \\ 
\texttt{[ 384, 5678 ]} &  & 1 & 16 \\ 
\texttt{[ 384, 5781 ]} &  & \textbf{4} & 22 \\ 
\texttt{[ 432, \ 520 ]} &  & 1 & 14 \\
\texttt{[ 432, \ 523 ]} &  & 1 & 20 \\
\texttt{[ 432, \ 734 ]} & \(\Hol(\ZZ_3 \times \ZZ_3)\) & 1 & 11 \\
\texttt{[ 480, 1189 ]} & \(S_4 \times \Hol(\ZZ_5)\) & 1 & 25 \\ 
\texttt{[ 486, \ \ 31 ]} & \(\Hol(\ZZ_{27})\) & 1 & 31 \\ 
\texttt{[ 486, \ 132 ]} & & \textbf{3} & 34 \\ 
\texttt{[ 500, \ \ 18 ]} & \(\Hol(\ZZ_{25})\) & 1 & 26 \\ 
\texttt{[ 506, \ \ \ 1 ]} & \(\Hol(\ZZ_{23})\) & 1 & 23 \\ 
	\end{tabular}
\end{table}

We also determined which finite simple groups have trivial Reidemeister spectrum. The Feit-Seitz theorem \cite[Thm.~C]{fs89-a} says that \(\Outc(G) = 1\) for a simple group \(G\), hence \(G\) has trivial Reidemeister spectrum if and only if \(\Out(G) = 1\). The finite simple groups that satisfy this are listed in \zcref{tbl:spectraSimple}; the data on the class numbers was taken from \cite{lueb00-a,abl06-a}.

\begin{table}
	\centering
	\small
	\caption{Simple groups \(G\) with trivial \(\Out(G)\)}
	\label{tbl:spectraSimple}
	\begin{tabular}{c||c|c|c}
		Class & Group \(G\)& Conditions &\(k(G)\) \\
		\hline
		\hline
		Symplectic groups& \(\mathrm{S}_{2n}(2)\) &\(n > 2 \)& see \cite[Thm.~3.7.3]{wall63-a}\\
		\hline
		& \(\mathrm{E}_7(2)\) &&\(531\) \\
		& \(\mathrm{E}_8(2)\) && \(1\,156\)\\
		& \(\mathrm{E}_8(3)\) &&\(12\,825\) \\
		\multirow{2}{*}{Exceptional groups}& \(\mathrm{E}_8(5)\) &&\(519\,071\) \\
		\multirow{2}{*}{of Lie type} & \multirow{2}{*}{\(\mathrm{E}_8(p)\)}&\multirow{2}{*}{\(p > 5\)}& \(p^8+p^7+2p^6+3p^5+10p^4\)\\
		&  && \(+16p^3+40p^2+67p+112\)\\
		& \(\mathrm{F}_4(3)\)& & \(273\)\\
		& \(\mathrm{F}_4(p)\)&\(p > 3\) & \(p^4+2p^3+7p^2+15p+31\)\\
		& \(\mathrm{G}_2(p)\)&\(p > 3\)&\(p^2+2p+9\) \\
		\hline
		\multirow{3}{*}{Mathieu groups} & \(\mathrm{M}_{11}\) && \(10\)\\
		& \(\mathrm{M}_{23}\) && \(17\)\\
		& \(\mathrm{M}_{24}\) && \(26\)\\
		\hline
		\multirow{2}{*}{Janko groups} & \(\mathrm{J}_{1}\) && \(15\)\\
		& \(\mathrm{J}_{4}\) && \(62\)\\
		\hline
		 & \(\mathrm{Co}_{1}\) && \(101\)\\
		Conway groups& \(\mathrm{Co}_{2}\) && \(60\)\\
		& \(\mathrm{Co}_{3}\) && \(42\)\\
		\hline
		Fisher group & \(\mathrm{Fi}_{23}\) && \(98\)\\
		\hline
		Rudvalis group& \(\mathrm{Ru}\) && \(36\)\\
		\hline
		Lyons group & \(\mathrm{Ly}\) && \(53\)\\
		\hline
		Thompson group & \(\mathrm{Th}\) &&\(48\)\\
		\hline
		Baby monster group & \(\mathrm{B}\) && \(184\) \\
		\hline
		Monster group & \(\mathrm{M}\) && \(194\)\\
		
	\end{tabular}
\end{table}	

Every group with \(\Autc(G) = \Aut(G)\) mentioned above has \(\Out(G) = 1\) (e.g. simple groups), \(Z(G) \neq 1\) (e.g. \(p\)-groups), or both (e.g. \SmallGroup{336}{210}). The existence of a group that does not fit in any of these categories is still an open question.
\begin{question}
	Does there exist a finite group \(G\) with trivial Reidemeister spectrum, \(\Out(G) \neq 1\) and \(Z(G) = 1\)?
\end{question}

\section{Trivial extended Reidemeister spectrum}\label{sec:trivespec}
We can extend the idea of having trivial Reidemeister spectrum to endomorphisms, taking into account that every group admits the trivial endomorphism \(g \mapsto 1\).
\begin{definition}
	A finite group \(G\) is said to have \emph{trivial extended Reidemeister spectrum} if
	\begin{equation*}
		\ESpecR{G} = \{1,k(G)\}.
	\end{equation*}
\end{definition}
Again, we immediately obtain a characterisation of the groups admitting this property.
\begin{proposition}
	A finite group \(G\) has trivial extended Reidemeister spectrum if and only if every endomorphism is either class-preserving or fixed-point-free.
\end{proposition}

\begin{question}
	Which finite groups have trivial extended Reidemeister spectrum?
\end{question}
Examples of groups with trivial extended Reidemeister spectrum include the cyclic groups \(\ZZ_p\) of prime order (no other examples of order \(< 512\) exist), as well as each of the non-abelian finite simple groups from \zcref{tbl:spectraSimple}. One may suspect from this that only simple groups have trivial extended Reidemeister spectrum. The \emph{quasisimple} groups, however, demonstrate that this is not the case.

\begin{definition}
	A group \(G\) is called \emph{quasisimple} if \(G\) is perfect, i.e. \(G = [G,G]\), and \(G/Z(G)\) is simple.
\end{definition}
We list some properties of a quasisimple group \(G\), obtained from \cite[Prop.~3.5]{tsan19-a}.
\begin{proposition}
	\label{prop:quasisimpleprops}
	Let \(G\) be a quasisimple group. Then:
	\begin{enumerate}[label=(\alph*)]
		\item A proper normal subgroup of \(G\) is contained in \(Z(G)\);
		\item The natural homomorphism \(\Aut(G) \to \Aut(G/Z(G))\) is injective;
		\item An endomorphism of \(G\) is either trivial or an automorphism.
	\end{enumerate}
\end{proposition}

The next \zcref[noref,nocap]{prop:quasisimplespectra} follows almost immediately from these properties.
\begin{proposition}
	\label{prop:quasisimplespectra}
	Let \(G\) be a quasisimple group such that \(G/Z(G)\) has trivial outer automorphism group. Then \(G\) has trivial (extended) Reidemeister spectrum.
\end{proposition}
Such quasisimple groups (that are not simple) exist, e.g. the groups 2.Co\textsubscript{1} (\(k(G) = 167\)), 2.Ru (\(k(G) = 61\)) and 2.B (\(k(G) = 247\)). The converse to \zcref{prop:quasisimplespectra} is not true: the quasisimple group 2.Sz(\(8\)) has trivial outer automorphism group and Reidemeister spectrum \(\{19\}\), however, its simple quotient Sz(\(8\)) has non-trivial outer automorphism group and Reidemeister spectrum \(\{5,11\}\). All groups in this paragraph were named using the notation from the ATLAS \cite{abl06-a}.

\medskip

Note that there are groups with trivial Reidemeister spectrum but non-trivial extended Reidemeister spectrum (e.g. \(S_n\), \(n \neq 6\)) and vice versa (e.g. \(\ZZ_p\), \(p \geq 3\)). The former family is expected to have many examples, but groups belonging to the latter family seem to be exceedingly rare. This raises the next question.

\begin{question}
	Let \(G\) be a finite group with trivial extended Reidemeister spectrum, but non-trivial Reidemeister spectrum, i.e. \(\SpecR{G} = \ESpecR{G} = \{1,k(G)\}\). Must \(G\) be isomorphic to \(\ZZ_p\) for some prime \(p \geq 3\)?
\end{question}

\section{Full (extended) Reidemeister spectrum}

The original definition of full Reidemeister spectrum is the following:
\begin{definition}
	A group \(G\) has \emph{full Reidemeister spectrum} if
	\begin{equation*}
		\SpecR{G} = \NN \cup \{\infty\}.
	\end{equation*}
\end{definition}
Examples of infinite groups with full Reidemeister spectrum include the free abelian groups of rank \(\geq 2\) \cite[Sec.~3]{roma11-a} and the free nilpotent groups of rank \(\geq 4\) and class \(2\) \cite[Sec.~4]{dtv20-a}. We extend this definition to finite groups by, once again, using that \(k(G)\) is an upper bound for Reidemeister numbers.
\begin{definition}
	A finite group \(G\) has \emph{full Reidemeister spectrum} if
	\begin{equation*}
		\SpecR{G} = \{1,\ldots, k(G)\}.
	\end{equation*}
\end{definition}
Infinite groups with full Reidemeister spectrum seem to be rather hard to find. For finite groups, one doesn't even need to search (see also \cite[Prop.~A.2.6]{send23-a}).
\begin{proposition}
	\label{prop:nofullreidspec}
	A non-trivial finite group \(G\) cannot have full Reidemeister spectrum.
\end{proposition}
\begin{proof}
	Let \(\varphi \in \Aut(G)\), which induces a permutation \(\Phi\) on \(\mathcal{C}(G)\). Clearly, a permutation on a set of \(n\) elements cannot fix exactly \(n-1\) elements. This means that \(\Phi\) cannot fix exactly \(k(G)-1\) conjugacy classes, and hence \(\varphi\) cannot have Reidemeister number \(k(G)-1\).
\end{proof}
We quickly move on to the extended Reidemeister spectrum.
\begin{definition}
	A group \(G\) has \emph{full extended Reidemeister spectrum} if
	\begin{equation*}
		\ESpecR{G} = \NN \cup \{\infty\}.
	\end{equation*}
\end{definition}
Examples include all finitely generated torsion-free nilpotent groups, as shown in \cite[Section 6]{dtv20-a}. We extend the definition of full extended Reidemeister spectrum to finite groups in the usual way.
\begin{definition}
	A finite group \(G\) has \emph{full extended Reidemeister spectrum} if
	\begin{equation*}
		\ESpecR{G} = \{ 1, \ldots, k(G)\}.
	\end{equation*}
\end{definition}
Again, we wonder which groups admit this property.
\begin{question}
	Which finite groups have full extended Reidemeister spectrum?
\end{question}

We shall start by excluding some families of groups from contention. Inspired by \zcref{prop:nofullreidspec}, we start by looking at Reidemeister number \(k(G)-1\).
\begin{lemma}
	\label{thm:kGminus1}
	Let \(G\) be a finite group of order \(> 2\) such that every non-trivial normal subgroup \(N\) of \(G\) intersects the centre \(Z(G)\) non-trivially. Then \(k(G)-1 \notin \ESpecR{G}\).
\end{lemma}

\begin{proof}
	We will work by contradiction. Let \(\varphi\) be an endomorphism of \(G\) with \(R(\varphi) = k(G) - 1\). Due to (the proof of) \zcref{prop:nofullreidspec}, \(\varphi\) cannot be an automorphism. Hence, the subgroup \(N_\varphi\) (as defined in \zcref{prop:goodnormalsubgroup}) is non-trivial and \(k(G/N_\varphi) \leq k(G)-1\). Let \(\bvarphi\) be the induced automorphism on \(G/N_\varphi\), then we find that 
	\begin{equation*}
		k(G)-1 = R(\varphi) = R(\bvarphi) \leq k(G/N_\varphi) \leq k(G)-1,
	\end{equation*}
	so all inequalities are actually equalities and \(k(G/N_\varphi) = k(G) - 1\). Only a single conjugacy class is lost when projecting from \(G\) to \(G/N_\varphi\), hence \(N_\varphi\) must be the union of exactly two conjugacy classes, one of which is the class of the identity. Then \(N_\varphi\) is a minimal normal subgroup, and since \(\ker(\varphi)\) is non-trivial and contained in \(N_\varphi\), these two subgroups are equal. We set \(N \assign N_\varphi = \ker(\varphi)\).
	
	Pick a non-trivial \(c \in N \cap Z(G)\). Since \(c \in Z(G)\), \([c]\) is a conjugacy class containing only \(c\); and since \(c \in N = \ker(\varphi)\), we have that \([\varphi(c)] = [1]\). Combining the assumption \(R(\varphi) = k(G) - 1\) with \zcref{thm:reidnrfixedpoints} tells us that there is precisely one conjugacy class with \([\varphi(g)] \neq [g]\); this conjugacy class must then be \([c]\). If \(g \notin \{1,c\}\), then \([\varphi(g)] = [g] \neq [1]\) and hence \(\varphi(g) \neq 1\). So \(N = \{1,c\}\).
	
	Suppose that \(\card{Z(G)} > 2\), then some \(x \notin \{1,c\}\) is contained in \(Z(G)\). Since \([\varphi(g)] = [g]\) for all \(g \neq c\), we get
	\begin{equation*}
		[cx] = [\varphi(cx)] = [\varphi(c)\varphi(x)] = [\varphi(x)] = [x].
	\end{equation*}
	But \(x\) is central, so its conjugacy class contains only \(x\) itself, and hence \(cx = x\). This implies that \(c = 1\), which is false, thus no such \(x\) can exist. We conclude that \(Z(G) = \{1,c\} = N\).
	
	By the first isomorphism theorem, \(G/N \cong \im(\varphi)\). Since \(\card{N} = 2\), \(\im(\varphi)\) is an index 2 subgroup of \(G\) and is therefore normal. Thus, \(\im(\varphi)\) intersects \(Z(G) = \{1,c\}\) non-trivially and hence \(c \in \im(\varphi)\). Let \(g \in G\) such that \(\varphi(g) = c\), then \(g \notin \ker(\varphi) = \{1,c\}\) and \([\varphi(g)] = [c] \neq [g]\), which contradicts what we deduced earlier.
\end{proof}
	This \zcref[noref,nocap]{thm:kGminus1} is applicable both to finite nilpotent groups (see \cite[\nopp{5.2.1.}]{robi96-a}) and quasisimple groups (see \zcref{prop:quasisimpleprops}). And, from \zcref{cor:oddordergroupreidnrodd}, it follows immediately that non-trivial finite groups of odd order cannot have full extended Reidemeister spectrum. All in all, we have eliminated the following families of groups.
	
\begin{theorem}
	Let \(G\) be a finite group of order \(>2\) that is either 
	\begin{itemize}
		\item of odd order,
		\item nilpotent, or,
		\item quasisimple,
	\end{itemize}
	then \(G\) cannot have full extended Reidemeister spectrum.
\end{theorem}

We have searched all groups of order \(< 1536\) and all groups with class number \(< 15\) for full extended Reidemeister spectra. Five such groups were found, they are listed in \zcref{tbl:fullextreidspec}.

\begin{table}[H]
	\centering
	\small
	\caption{Known finite groups with full extended Reidemeister spectrum}
	\label{tbl:fullextreidspec}
	\begin{tabular}{c|c|c}
		SmallGroup ID & Structure &  \(k(G)\) \\\hline
		\SmallGroup{\ \ 1}{\ \ 1} & \(1\) & \(1\) \\
		\SmallGroup{\ \ 2}{\ \ 1} & \(\ZZ_2\)  & \(2\) \\
		\SmallGroup{\ \ 6}{\ \ 1} & \(S_3\)  & \(3\) \\
		\SmallGroup{\ 12}{\ \ 3} & \(A_4\)  & \(4\) \\
		\SmallGroup{\ 72}{\ 41} & \(M_9\) & \(6\)
	\end{tabular}
\end{table}

Given that no groups of order \(72 < \card{G} < 1536\) or class number \(6 < k(G) < 15\) have full extended Reidemeister spectrum, we suspect the groups in \zcref{tbl:fullextreidspec} are the only ones with this property.
\begin{conjecture}
	The groups in \zcref{tbl:fullextreidspec} are the only finite groups with full extended Reidemeister spectrum.
\end{conjecture}

\section*{Acknowledgement}
The author would like to thank the anonymous referee for their careful reading and their helpful remarks and suggestions.

\emergencystretch=1em
\printbibliography

\end{document}